\documentclass[conference]{IEEEtran}
\IEEEoverridecommandlockouts

\usepackage{cite}
\usepackage{amsmath,amssymb,amsfonts}
\usepackage{algorithmic}
\usepackage{graphicx}
\usepackage[dvipsnames]{xcolor}
\usepackage[final]{microtype}
\usepackage[italic]{mathastext}
\usepackage{libertine}
\usepackage[T1]{fontenc}
\usepackage{textcomp}
\usepackage[varqu,varl]{zi4}
\usepackage[all]{nowidow}
\usepackage[keeplastbox]{flushend}
\usepackage{fancyhdr}

\usepackage{multirow}
\usepackage[hyphens]{url}
\usepackage{subfigure}
\usepackage{booktabs}

\usepackage{listings}

\definecolor{codegreen}{rgb}{0,0,0.6}
\definecolor{codegray}{rgb}{0.5,0.5,0.5}
\definecolor{codepurple}{rgb}{0.58,0,0.82}
\definecolor{backcolour}{rgb}{0.95,0.95,0.92}

\lstdefinestyle{lststyle}{
  backgroundcolor=\color{backcolour}, commentstyle=\color{codegreen},
  keywordstyle=\color{magenta},
  numberstyle=\tiny\color{codegray},
  stringstyle=\color{codepurple},
  basicstyle=\ttfamily\footnotesize,
  breakatwhitespace=false,         
  breaklines=true,                 
  captionpos=b,                    
  keepspaces=true,                 
  numbers=left,                    
  numbersep=5pt,                  
  showspaces=false,                
  showstringspaces=false,
  showtabs=false,                  
  tabsize=2
}

\def\BibTeX{{\rm B\kern-.05em{\sc i\kern-.025em b}\kern-.08em
    T\kern-.1667em\lower.7ex\hbox{E}\kern-.125emX}}

\begin{document}


\title{Design and accuracy trade-offs in Computational Statistics}

\author{\IEEEauthorblockN{Tiancheng Xu\IEEEauthorrefmark{1}\thanks{* Currently with Google (tcxxxx@google.com). This work was completed during his PhD study at Rice University.}}
\IEEEauthorblockA{
\textit{Rice University}\\
Houston, USA \\
txu@rice.edu}
\and
\IEEEauthorblockN{Alan L. Cox}
\IEEEauthorblockA{\textit{Rice University} \\
Houston, USA \\
alc@rice.edu}
\and
\IEEEauthorblockN{Scott Rixner}
\IEEEauthorblockA{\textit{Rice University} \\
Houston, USA \\
rixner@rice.edu}
}

\maketitle
\pagestyle{plain}


\begin{abstract}

Statistical computations are becoming increasingly important. These computations often need to be performed in log-space because probabilities become extremely small due to repeated multiplications. While using logarithms effectively prevents numerical underflow, this paper shows that its cost is high in performance, resource utilization, and, notably, numerical accuracy. This paper then argues that using posit, a recently proposed floating-point format, is a better strategy for statistical computations operating on extremely small numbers because of its unique encoding mechanism. To that end, this paper performs a comprehensive analysis comparing posit, binary64, and logarithm representations, examining both individual arithmetic operations, statistical bioinformatics applications, and their accelerators. FPGA implementation results highlight that posit-based accelerators can achieve up to two orders of magnitude higher accuracy, up to 60\% lower resource utilization, and up to $1.3\times$ speedup, compared to log-space accelerators. Such improvement translates to $2\times$ performance per unit resource on the FPGA.

\end{abstract}

\begin{IEEEkeywords}
Floating-point arithmetic, numerical analysis, Posit, FPGA acceleration, statistical bioinformatics, computational statistics
\end{IEEEkeywords}

\section{Introduction}
\label{sec:intro}

Many modern applications, such as bioinformatics and finance, increasingly rely on statistical computations. Generally, such statistical computations iteratively update probabilistic states, where the current state depends on prior states.  These iterative updates involve the addition and multiplication of probabilities. The multiplications decrease the state probabilities with each successive iteration. Two classic examples of computations that follow this pattern are the construction of Hidden Markov Models (HMM)~\cite{hmm_tutorial} and the Poisson Binomial Distribution~\cite{tang2023poisson, hong2013computing}.

Due to the iterative computation pattern, probabilities in such applications can easily become extremely small and thus underflow. The 64-bit IEEE floating-point (binary64) has insufficient dynamic range: the smallest positive representable binary64 is only $2^{-1,074}$. This is not small enough in many situations. For example, using the forward algorithm to build an HMM on $500,000$ site long Human-Chimp-Gorilla (HCG) genome sequences~\cite{hcg_data} can yield likelihoods as low as $2^{-2,900,000}$~\cite{liu2021variational}. Nonetheless, such likelihoods must be preserved because underflow to zero prevents proper convergence and leads to incorrect results in algorithms such as Variational Inference and Markov Chain Monte Carlo~\cite{liu2021variational, blei2017variational,10.1093/sysbio/syx085}.

In order to overcome the range limitations of IEEE floating-point numbers, statistical computations are often performed in log-space~\cite{carpenter2017stan,hmm_numerical, liu2021variational,n_gram_models, durbin1998biological, eddy2011accelerated, stanford_cs109,blei2017variational, wilm2012lofreq, rosenman2019some, 10.1093/sysbio/syx085,10.1093/bioinformatics/btm388}. Using the logarithm of the probability, instead of the probability itself, converts extremely small numbers into negative numbers that are well within binary64's representable range. For example, the natural logarithm of $2^{-2,900,000}$ is $-2,010,126.824$, which can be easily represented with binary64. Therefore, the use of logarithms is the standard approach to perform statistical computations.

However, we find that this actually creates a trade-off, because numerical accuracy is a function of both the precision and the range of the numerical format.  While log-space clearly increases the range of representable numbers, this paper shows that it does so at the cost of reduced precision and numerical accuracy.

For example, this paper reveals that, compared to binary64, using logarithms only improves accuracy \textbf{outside} binary64's normal range (Figure~\ref{fig:ops_boxplot}). Within binary64's normal range, using logarithms actually results in \textbf{worse} numerical accuracy. This is because converting a probability number into log-space wastes available exponent bits. The logarithm value's exponent typically requires far fewer bits to encode, and thus, the exponent bits go unused. Meanwhile, the fraction bits in the logarithm value are effectively used to encode both the fraction and the exponent of the original value, reducing the number of bits for precision.


This paper highlights that the recently proposed floating-point format, posit~\cite{gustafson2017beating}, can achieve both high numerical accuracy and wide dynamic range simultaneously. Instead of having fixed bit fields, posit dynamically adjusts the number of exponent and fraction bits based on needs. When numbers become extremely small, posit devotes more bits to the exponent so that such numbers can be represented without underflow. When fewer exponent bits are needed, posit uses more bits as fraction bits to achieve higher precision, and thus, higher numerical accuracy. This paper shows that posit achieves higher numerical accuracy than logarithms on numbers where binary64 underflows to zero.


Besides numerical accuracy, using logarithms also hurts performance and hardware resource costs. In log-space, while a multiplication becomes slightly simpler, an addition becomes a sequence of logarithm and exponential operations. Such complications not only increase critical path latency, but also require the implementation of more expensive logarithm and exponential operators. Evaluation in this paper shows that, compared to binary64, log-space addition is 10$\times$ slower and requires 8$\times$ as many LUTs and FFs on an FPGA. Meanwhile, this paper demonstrates that posit shows great advantages in performance and resource cost by avoiding the need to operate in log-space.

Built upon our analysis of the number formats at the arithmetic level, this paper further examines using posit in two critical bioinformatics applications operating on small probabilities: VICAR and LoFreq. VICAR is a phylogenetics application that uses HMMs to analyze evolutionary parameters of species trees. LoFreq is a genomics tool that identifies genome variants using the Poisson Binomial Distribution to analyze genome alignment data. The state-of-the-art software implementations of both applications suffer from long execution times due to the use of logarithms. This paper builds highly optimized FPGA accelerators for both applications and then studies the application-level numerical accuracy, performance, and hardware resource cost of the accelerators. Our evaluation results highlight that using posit leads to improvements in all metrics.

In particular, using posit leads to two orders of magnitude higher numerical accuracy in final application results. Besides, posit-based accelerators can achieve up to $33\%$ higher performance and $60\%$ lower resource utilization, compared to highly parallelized log-space accelerators. This results in gains of up to a factor of 2 in performance per resource unit on the FPGA.
In summary, the main contributions of this paper include:
\begin{itemize}
\item This paper is the first to perform an in-depth analysis comparing the numerical accuracy of binary64, logarithm, and posit together. The key observations include, for example, that posit is more accurate than using logarithms outside binary64's range, and that the numerical accuracy of posit changes more steadily than logarithms.

\item This paper builds highly optimized FPGA accelerators using posit for critical statistical applications, and demonstrates the benefits of using posit.

\item This paper, for the first time, provides the insight that posit is well-suited for statistical computations and, thus, should be used in future architectures for such applications.

\end{itemize}


The rest of this paper proceeds as follows: Section~\ref{sec:problem} describes the numerical underflow challenge in statistical computations and the standard log-based solution. Section~\ref{sec:posit} introduces posit, and Section~\ref{sec:tradeoff} presents an in-depth numerical analysis of all number formats. Section~\ref{sec:case} introduces real applications and hardware accelerators for a case study, and Section~\ref{sec:eval} highlights the benefits of using posit. Section~\ref{sec:related} discusses related works and Section~\ref{sec:conc} concludes the paper.

\section{Problem}
\label{sec:problem}

\subsection{Motivation}
Statistical computations are increasingly common in modern applications, where fundamental operations are the multiplication and addition of probabilities. In these applications, a key challenge is that numbers tend to become extremely small, and such small values must be preserved. As probabilities are numbers between $0$ and $1$, repeated multiplications can make them decrease quickly.

The IEEE floating-point standard (IEEE 754) defines the format of double-precision floating-point numbers (binary64)~\cite{ieee754_2019}, which is the number type with the highest precision and largest dynamic range that is common to all current machine architectures. The smallest positive number binary64 can represent is $2^{-1,074}$, which is not small enough in many cases. Consider the binomial distribution as an example, to compute the probability $P$ of observing $N$ successes in $N$ Bernoulli trials where the success rate is $0.3$ ($P = 0.3^N$), $P$ underflows for any $N$ larger than $618$ if represented in binary64.

Statistical computations often have an iterative pattern where current states rely on prior ones. Multiply is commonly part of the iterations. Two classic examples are the Hidden Markov Model (HMM) and Poisson Binomial Distribution (PBD). In the HMM forward algorithm, the core computation can be expressed as:

\[alpha\_i = alpha\_i^{\prime} \times a \times e\]

\noindent In this equation, all variables are probabilities.
$alpha\_i$ and $alpha\_i^{\prime}$ denote the $i^{th}$ $alpha$ state in the current ($t$) and previous ($t-1$) iterations, respectively. In computing the probability mass function of a PBD, the core computation can be expressed as:

\[p_{k} = {p_{k}}^{\prime} \times (1.0 - p_{n}) + {p_{k-1}}^{\prime} \times p_{n}\]

\noindent In this equation, all operands are probabilities. $p_{k}$ denotes the $k^{th}$ state in the current iteration, and ${p_{k}}^{\prime}$ and ${p_{k-1}}^{\prime}$ denote states from the previous iteration. Both algorithms exhibit the two traits: iterative computation and repeated multiplications.

\begin{figure}[t]
    \centering
    \includegraphics[width=0.95\columnwidth]{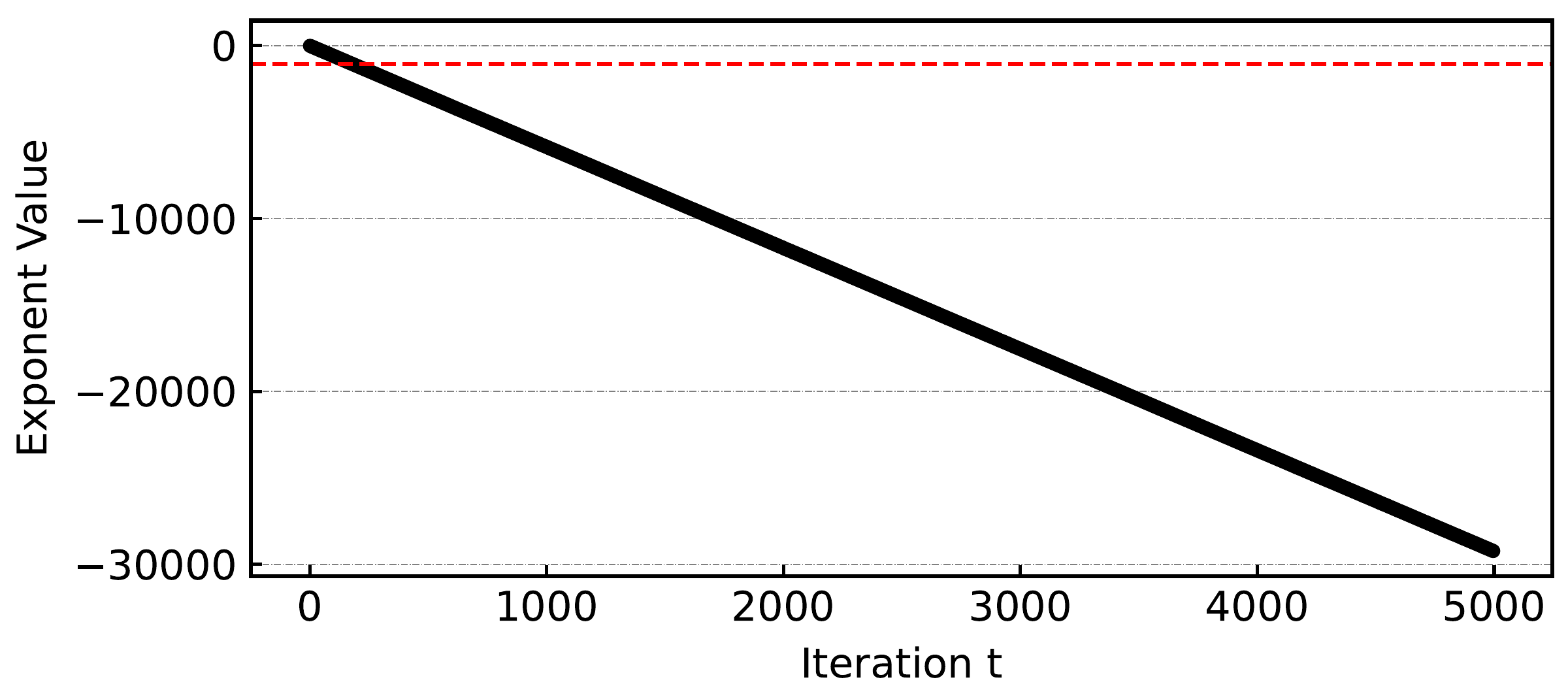}
    \caption{Base-2 exponent value of $alpha$ over iterations.}
    \label{fig:exp_over_time}
    \end{figure}

Figure~\ref{fig:exp_over_time} shows how $alpha$ changes over the course of running the forward algorithm to illustrate the scale of numbers in such computations. The X axis shows the number of iterations, and the Y axis shows the base-2 exponent of $alpha$. The experiment is done using the MPFR arbitrary precision library so that the exact exponent can be tracked even when numbers become extremely small. The dotted red line shows the exponent of the smallest positive binary64. If the computation is implemented using binary64, most results would underflow and become absolutely meaningless.

\subsection{Current Solution: Using Logarithms}
\label{sec:statusquo}
In many statistical applications, such extremely small probabilities are significant, so their numerical values must be preserved. In such scenarios, the standard practice is to perform the computation in log-space using log probabilities or log likelihoods~\cite{stan_doc_log,
carpenter2017stan, hmm_numerical, liu2021variational,
n_gram_models, durbin1998biological, eddy2011accelerated, stanford_cs109,
blei2017variational, wilm2012lofreq, rosenman2019some, 10.1093/sysbio/syx085,
10.1093/bioinformatics/btm388}.
In the probabilistic programming language \textit{Stan}, all probability computations operate in log-space~\cite{carpenter2017stan, stan_doc_log}. Using logarithms effectively prevents underflow by significantly expanding the range of representable numbers without using more bits. When converted to log-space, an extremely small probability number between 0 and 1 in the original linear space becomes a normal negative number. For example, the log of an extremely small number $2 ^ {-120,000}$ is approximately $-83177.66$, which can easily be represented using an IEEE floating-point number. With logarithms, the smallest positive number that can be represented in binary64 becomes approximately $2^{-2.59 \times 10 ^ {308}}$ instead of $2^{-1,074}$. This expands the dynamic range to effectively infinite, making numerical underflow almost impossible.

However, the increased dynamic range comes at the cost of more complex computations. While multiply becomes add, add becomes a series of expensive operations. For example, consider two numbers, $x$ and $y$. In log-space, these numbers become $lx$ and $ly$, where $lx = log(x)$ and $ly = log(y)$. The simple addition of these numbers, $x + y$, now becomes:

\begin{equation}
    log(x+y)=log(exp(lx) + exp(ly))
    \label{eq:lse_naive}
\end{equation}
    
Furthermore, the computation is even more complex than equation~\ref{eq:lse_naive} as described below, because directly evaluating the exponential terms often causes numerical underflow and overflow issues. Assuming all numbers are represented in binary64, exponentiating a value smaller than $-745.133$ will underflow, and exponentiating a value greater than $709.782$ will overflow. Therefore, naively adding $x$ and $y$ in log-space is numerically unstable.

The \textbf{Log Sum Exp} (LSE) technique is generally used instead to avoid such numerical issues~\cite{compute_log_sum_exp,lse_blog,cs231n_softmax,blei2017variational}. The intuition is to make the input to the exponential operations closer to 0. The mathematical form of LSE is as follows:

\begin{equation}
    \begin{split}
    &m = max(lx, ly) \\ 
    &log\_sum\_exp(x, y) = \\
    &\quad \quad \quad m + log(exp(lx - m) + exp(ly - m))
    \end{split}
    \label{eq:lse_two}
\end{equation}

While Equations~\eqref{eq:lse_naive} and~\eqref{eq:lse_two} are mathematically equivalent, the LSE technique completely eliminates overflow and greatly reduces the chance of underflow. With the subtraction, the input to $exp$ will always be less than or equal to $0$; thus, overflow will never happen. Underflow can also be avoided as long as $lx$ and $ly$ are relatively close. For example, if $lx=-1,000$ and $ly=-999$, both $exp(lx)$ and $exp(ly)$ will underflow in Equation~\eqref{eq:lse_naive}; meanwhile, Equation~\eqref{eq:lse_two} can compute the correct result without underflow.

More generally, LSE works for calculating the sum of multiple numbers represented in log-space. For $s = {x_{1}, ... , x_{N}}$ represented by their log values $ls = {lx_{1}, ... , lx_{N}}$, the sum operation $\sum_{n=1}^{N}x_{n}$ using LSE is shown in Equation~\eqref{eq:lse}:

\begin{equation}
    \begin{split}
    &m = max(lx_{1}, ... lx_{N}) \\
    &log\_sum\_exp(x_{1}, ..., x_{N})= m + log\sum_{n=1}^{N}exp(lx_{n} - m)
    \end{split}
\label{eq:lse}
\end{equation}

\subsection{Drawbacks of Using Logarithms}
\label{sec:log_drawback}
Using logarithms has been the standard approach because it enables representing extremely small numbers and is easy to implement in software. However, these benefits come at the cost of performance, hardware cost, and numerical accuracy.

The LSE operation not only takes many more cycles compared to addition, but also prevents straightforward parallelization. Besides the expensive logarithm and exponential operations, the max operation introduces another synchronization point in the dataflow.

Surprisingly, we find that numerical accuracy is hurt by using logarithms. This is because using logarithms reduces precision. Recall that precision is dictated by the number of fraction bits. As numbers operate in log-space, fraction bits are used to encode both the fraction and the exponent of the original value. Thus, fewer fraction bits are actually used for precision. From another perspective, most logarithm values in these computations have relatively small positive exponents, e.g., $8$. As a result, most exponent bits are wasted.

To illustrate, we use two logarithm values, $-402.1$ and $-408.1$ (encoded in binary64), as an example. Both have the same exponent. However, their original values are $1.856\times2^{-581}$ and $1.178\times2^{-589}$, respectively, where the exponent difference is actually large. This shows that the order of magnitude information (exponent) is effectively encoded in the fraction bits when using logarithms. From another perspective, the exponent binary bits of the log value $-402.1$ and its original value $1.856\times2^{-581}$ are $10000000111$ and $00110111010$, respectively. When using logarithms within binary64's range, most of the exponent bits in the log values are unused as $0$s, and thus, wasted. Note that the exponent of the log value is only $8$. The MSB in its exponent bits is used only because of the bias ($-1023$).

Our quantitative analysis in section~\ref{sec:tradeoff_acc} will show that arithmetic operations become less accurate when they are done in log-space compared to binary64.
\section{Posit}
\label{sec:posit}

\begin{figure}[t]
  \subfigure[IEEE Floating-Point: each field has a fixed length.]
  {
      \includegraphics[width=0.4\textwidth]{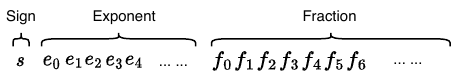}
      \centering
      \label{fig:ieee_bits}
  }
  \subfigure[Posit: all fields except for the sign have a variable length.]
  {
      \includegraphics[width=0.4\textwidth]{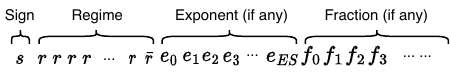}
      \centering
      \label{fig:posit_bits}
  }
  \caption{IEEE Floating-Point and Posit Number Formats.}
  \label{fig:number_formats}
\end{figure}

This paper argues that \textit{posit}~\cite{gustafson2017beating}, a recently proposed floating-point number format, is well suited for statistical computations operating on extremely small probabilities. Posit achieves both wide dynamic range and high numerical accuracy thanks to its number encoding mechanism that is fundamentally different from the IEEE standard.

The IEEE floating-point format is suboptimal for statistical computations. Figure~\ref{fig:ieee_bits} shows the format of an IEEE floating-point number. All bit fields have a \textbf{fixed} length. For example, a binary64 has one sign bit, $11$ exponent bits, and $52$ fraction bits. This is painful when more exponent bits are needed, such as when representing $2^{-10,000}$. In such cases, IEEE numbers must be converted to log-space to prevent underflow.

In contrast, posit avoids such problems by design. Posit dynamically adjusts the allocation of bits to the different fields within the number based on needs. When fewer exponent bits are needed, the ``spare'' bits can be used as fraction bits, which leads to higher precision and eventually higher numerical accuracy. Meanwhile, more bits can be used as exponent bits on demand: this enables posit to have a significantly wider dynamic range, allowing it to represent much smaller numbers without using logarithms. Note that such a mechanism is achieved in hardware and is transparent to software.

\textbf{Posit Format}. A posit number has two configuration parameters: the total number of bits ($N$) and the maximum number of exponent bits ($ES$), denoted in this paper as posit$(N, ES)$. Figure~\ref{fig:posit_bits} shows the bit composition of an $N$-bit posit with $ES$ exponent width. The sign bit works the same way as IEEE floating-point. A notable difference is the \textbf{regime bits}. The regime bits, together with the exponent bits that immediately follow, serve the same function as the exponent bits in IEEE floating-point format. As shown in Equation~\ref{eq:posit_val}, the value $v$ of a posit number is computed in the following way, where $e$ is an unsigned integer value encoded by the exponent bits, $f$ is the fraction value encoded by the fraction bits, and $l$ is the number of times $r$ is repeated in the regime bits.

\begin{equation}
  \begin{split}
    &useed = 2 ^ { 2 ^ {ES} } \\
    &k = 
    \left\{
      \begin{aligned}
        & -l  && \text{if } r = 0 \\
        & l-1 && \text{if } r = 1 \\
      \end{aligned}
      \right. \\
    &v = 
      \left\{
          \begin{aligned}
            & \text{NaR}  && \text{if } p = 100...0 \\
            & 0  && \text{if } p = 000...0 \\
            & (-1) ^ {sign} \times useed ^ {k} \times 2 ^ {e} \times (1 + f) && \text{otherwise} \\
          \end{aligned}
          \right.
  \end{split}
  \label{eq:posit_val}
  \end{equation}

\begin{table}[t]
  \centering
  \caption{Dynamic Range and Precision of Different Number Formats.}
  \label{tab:num_es}
  \begin{tabular}{l | l l c }
  \toprule
   Format & \vtop{\hbox{\textit{useed}}} & \vtop{\hbox{\strut Smallest Representable}\hbox{\strut Positive Number}} & \vtop{\hbox{\strut Max Num. of}\hbox{\strut Fraction Bits}} \\ [0.5ex]
  \midrule
   binary64  & - & $2 ^ {-1,074}$ & 52 \\
   \midrule
   posit(64,6) & $2 ^ {64}$ & $2 ^ {-3,968}$ & 55 \\
   posit(64,9) & $2 ^ {512}$ & $2 ^ {-31,744}$ & 52 \\
   posit(64,12) & $2 ^ {4,096}$ & $2 ^ {-253,952}$ & 49 \\
   posit(64,15) & $2 ^ {32,768}$ & $2 ^ {-2,031,616}$ & 46 \\
   posit(64,18) & $2 ^ {262,144}$ & $2 ^ {-16,252,928}$ & 43 \\
   posit(64,21) & $2 ^ {2,097,152}$ & $2 ^ {-130,023,424}$ & 40 \\
   \bottomrule
  \end{tabular}
\end{table}

As these equations show, the $useed$ value and the regime bits combine to form an exponential scaling factor for the value, $v$, of the posit number. The $useed$ value is determined by the number of exponent bits, $ES$, in the posit configuration. Recall that the regime bits consist of a sequence of $l$ identical bits ($r$) ended by the opposite bit ($\bar{r}$). The $k$ value is decided by the number of times $r$ is repeated in the regime bits, $l$, within the particular posit number, as shown in the above equations. For example, if the regime bits are $0001$, then $k=-3$. The minimum number of regime bits is 2 ($01$ or $10$), and the maximum number is $N-1$ (leaving only 1 bit for sign). Thus, extreme posit numbers may have only regime bits and no exponent or fraction bits at all.

Exponent bits immediately follow the regime bits: the next $ES$ bits are exponent bits when there are at least $ES$ remaining; all remaining bits are exponent if there are fewer than $ES$ bits remaining. All remaining bits, if any, after the exponent bits are fraction bits.

Unlike IEEE floating-point, posits do not use a bias when decoding the exponent bits, and $e$ is an unsigned number based on the exponent bits. Furthermore, the implicit bit is always $1$ for fraction bits because there are no subnormal numbers. There are only two special posit values: $0$ and Not a Real (NaR). As shown in Equation~\eqref{eq:posit_val}, the bit pattern of all $0$s represents $0$ and $100...0$ ($1$ followed by all $0$s), represents NaR. Unlike IEEE floating-point where there are positive and negative zero, there is only one zero in posit. Both infinity and Not a Number in IEEE floating-point are represented by NaR in posit.

\textbf{Example}. Consider the posit$(8,2)$ bit sequence $0\_0001\_10\_1$ (underscores are added to make it easier to see the different fields). The sign bit is zero, so this is a positive number. The regime bits are $0001$, so $k=-3$. The exponent bits are $10$, so $e=2$. The remaining bit, $1$, is the only fraction bit, and thus, the significand value (fraction plus the implicit $1$) is $1.1_2$ ($1.5$). Following Equation~\eqref{eq:posit_val}, the final value is $1.5 \times 2 ^ {-10}$, computed from $((2 ^ {4}) ^ {-3} \times 2 ^ {2} \times 1.5)$. Note that $ES$ is $2$ in this example. When $ES$ changes, the decoded value would be different for the same bit sequence.

\textbf{Posit Configuration}. The number of exponent bits, $ES$, is a key configuration parameter that allows posits to have different dynamic ranges. Since this paper focuses on probabilities (numbers between 0 and 1), the smallest representable positive number is used to measure dynamic range. This value of a given posit configuration is computed from multiplying \textit{useed} by the minimum value of $k$, which is always $-62$ when $N=64$. Therefore, when $N$ is fixed, posits with a larger $ES$ always have a wider dynamic range. Table~\ref{tab:num_es} shows how the range expands as $ES$ increases.

$ES$ also impacts posits' precision, but in a more subtle way. A larger $ES$ reduces the max number of bits available for fraction, as shown in Table~\ref{tab:num_es}. However, choosing a larger $ES$ can also increase precision. This is because a larger $ES$ can reduce the number of regime bits used, saving more bits for the fraction. For example, to encode $2^{-2,048}$, posit$(64,6)$ needs $33$ regime bits ($k=-32$), leaving only $24$ bits for the fraction. Meanwhile, posit$(64,9)$ needs only $5$ regime bits, leaving $49$ bits for the fraction. Thus, a careful quantitative study of how $ES$ impacts precision and numerical accuracy is required.

This paper uses three posit configs, posit(64,9), posit(64,12), and posit(64,18), for the analysis in later sections. Each is chosen for a reason. Posit(64,9) was selected for a direct comparison to binary64: posit(64,9) offers up to 52 fraction bits, matching binary64's precision, while providing a much wider dynamic range. Posit(64,18) was selected because it has a sufficient range for extremely small numbers as observed in critical statistical bioinformatics applications. Posit(64,12) represents an option balancing precision and range.
\section{Quantitative Trade-off Analysis}
\label{sec:tradeoff}

This section provides a comprehensive analysis of fundamental arithmetic operations using posits, logarithms, and binary64. The bit-width of all compared formats is fixed to $64$ for a fair and direct comparison. The experimental results and insights presented here apply universally across applications and constitute a core contribution of this paper.

\begin{figure*}[ht]
    \subfigure[Addition.]
    {
        \centering
        \includegraphics[width=1.95\columnwidth]{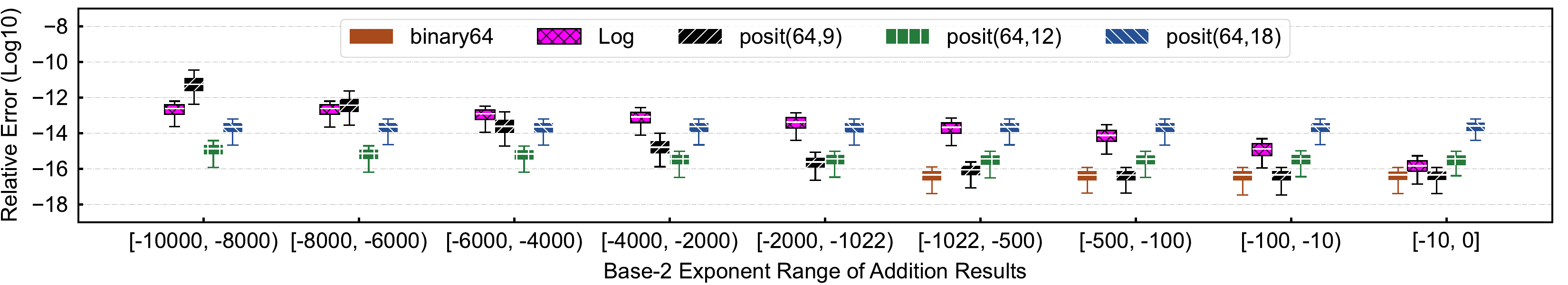}
        \label{fig:box_add}
    }
    \subfigure[Multiplication.]
    {
        \centering
        \includegraphics[width=1.95\columnwidth]{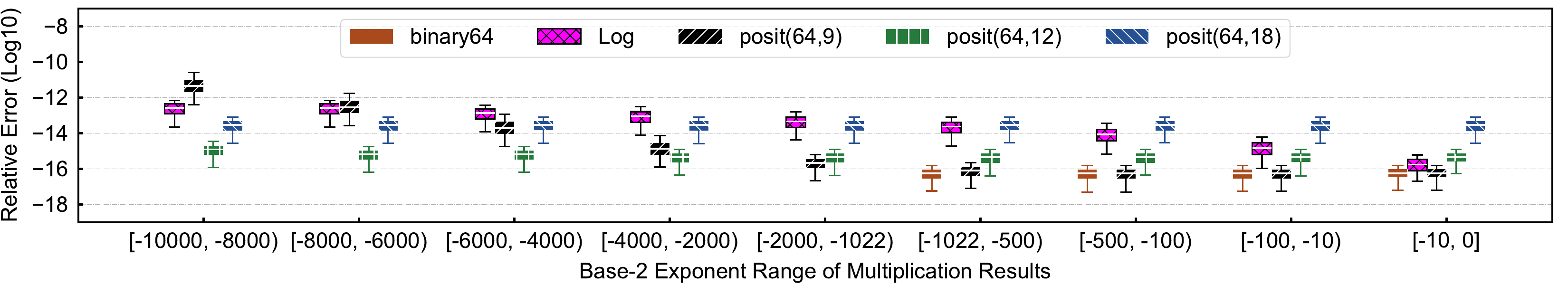}
        \label{fig:box_mul}
    }
    \caption{Individual operation accuracy on numbers of different magnitudes.}
    \label{fig:ops_boxplot}
\end{figure*}

\subsection{Numercial Accuracy}
\label{sec:tradeoff_acc}
This section compares the numerical accuracy of fundamental arithmetic operations using binary64, logarithms, and three posit configurations. Numerical accuracy measures how close a computed value is to the mathematically correct value. In this analysis, the MPFR library is used to compute the correct values~\cite{mpfr_paper,mpfr_website}. Results from 256-bit MPFR are regarded as the baseline correct values~\cite{debug_posit_pldi}.

The input operands are collected from both a real phylogenetics application and uniform sampling implemented in MPFR. Then, operands are converted from MPFR to each 64-bit format to perform the arithmetic operation. For logarithms, operands are transformed into log-space in MPFR. When the operation finishes, results from each format are converted back to MPFR to calculate accuracy. The relative error $|\frac{x-y}{x}|$ is computed to measure accuracy, where $x$ is the 256-bit MPFR result and $y$ is the 64-bit arithmetic result.

\noindent \textbf{Results.}
Figure~\ref{fig:ops_boxplot} presents the accuracy of individual add and multiply operations in different formats. The x-axis shows the exponent of the operands. This exponent corresponds to the exponent value in an IEEE floating-point number, but not to the value of $e$ in a posit number. Both figures are drawn from the same experimental data of $1,000,000$ add and $550,000$ multiply, where operation results range from $2^{-10,000}$ to $1$.

The y-axis shows the relative error on a log10 scale. Recall that accuracy is measured by relative errors. Each rectangular box represents one format. Overall, the accuracy of a format is higher when the box is lower. The rectangle's top and bottom lines are the $75^{th}$ and $25^{th}$ percentiles, respectively. The horizontal line within the rectangle shows the median. The whiskers are the $95^{th}$ and $5^{th}$ percentiles, respectively. Binary64 is not shown in ranges to the left of $2^{-1,022}$ due to underflow and having large errors in its subnormal range.

There are several key takeaways from the analysis. \textbf{First, using logarithms leads to worse numerical accuracy within binary64's normal range.} Within binary64's normal range (exponent $-1,022$ to $0$), logarithm's accuracy gets worse as numbers decrease. This confirms our discussion in Section~\ref{sec:log_drawback} that using logarithms reduces precision. Recall that accuracy is dictated by precision and dynamic range. Since these numbers are all within representable range when using logarithms, the loss of accuracy must be due to loss of precision, not range.

\textbf{Second, using posit leads to higher numerical accuracy than using logarithms.} Outside binary64's normal range, all three posits have higher accuracy, except for posit$(64,9)$ in the range of $[-10,000, -6,000)$, where it uses a large number of regime bits. Within binary64's normal range, posit$(64,9)$'s accuracy is constantly higher.

\textbf{Third, posits achieve the best of both worlds.} In ranges where binary64 underflows, posits maintain high numerical accuracy. Compared to logarithms, posits have higher overall accuracy. Compared to binary64, posits have high accuracy on a much wider range of numbers.

    \begin{table}[t]
      \centering
      \caption{Resource Utilization of Individual Arithmetic Units.}
      \resizebox{1\columnwidth}{!}
      {
      \label{tab:op_resource}
      \begin{tabular}{l | r r r r r}
        \toprule
       Arithmetic Unit & LUT & Register  & DSP & \vtop{\hbox{\strut Clock}\hbox{\strut Cycle}} & \vtop{\hbox{\strut Max Clock}\hbox{\strut Frequency (MHz)}} \\ [0.5ex]
       \midrule
       binary64 add & 679 & 587 & 0 & 6 & 480 \\
       Log add (binary64 LSE) & 5,076 & 5,287 & 34 & 64 & 346 \\
       posit(64,12) add & 1,064 & 1,005 & 0 & 8 & 354 \\
       posit(64,18) add & 1,012 & 974 & 0 & 8 & 358 \\
       \midrule
       binary64 mul & 213 & 484 & 6 & 8 & 480 \\
       Log mul (binary64 add) & 679 & 587 & 0 & 6 & 480 \\
       posit(64,12) mul & 618 & 1,004 & 9 & 12 & 336 \\
       posit(64,18) mul & 558 & 969 & 10 & 12 & 336 \\
       \bottomrule
      \end{tabular}
      }
    \end{table}

\subsection{Performance and Resource}
\label{sec:tradeoff_resource}
This section compares the latency, resource, and clock frequency on an FPGA of adder (add) and multiplier (mul) units of different formats. The comparison is done on an FPGA because software-emulated posit is too slow for practical use.

Note that the binary64 LSE operates on two inputs as in Equation~\eqref{eq:lse_two}. The experiment is performed on an Xilinx Alveo U250 FPGA. Binary64 add, binary64 mul, and LSE all use Xilinx's standard, optimized floating-point library, LogiCORE IP v7.1~\cite{logicoreip}. The posit units are implemented using a state-of-the-art implementation (MArTo)~\cite{uguen2019evaluating}. MArTo is a rigorously optimized HLS implementation that has shown better performance and usability over previous prototypes~\cite{podobas2018hardware,jaiswal2019pacogen,chaurasiya2018parameterized}. All units are placed and routed by Xilinx Vivado 2020.2. Table~\ref{tab:op_resource} shows the post-routing resource and latency of all units when they operate at max clock frequency. Conclusions persist when all units operate at the same $300$ MHz frequency.

First and foremost, posit adders consume fewer resources and have lower latency than binary64 LSE. This is because the log and exp operations in LSE are expensive in both resources and clock cycles. However, compared to binary64 adders, posit adders use more resources. More concretely, a posit$(64,12)$ adder consumes $70.3\%$ more LUTs and $44.0\%$ more registers than a binary64 adder. This is consistent with results from previous studies~\cite{uguen2019evaluating, murillo2023generating}. Remember, though, that MArTo is an HLS-based research prototype, whereas the LogiCORE IP is an optimized, RTL-based industrial product. Regardless, managing the variable-length fields in posits introduces extra overhead. MArTo also uses an internal data type that is larger than 64 bits for correct rounding~\cite{uguen2019evaluating}. Finally, compared to binary64, posits take an extra 2--4 cycles. This is due to the overhead of HLS and conversion between internal data types.

\section{Case Studies}
\label{sec:case}
The remainder of this paper studies how findings at the arithmetic op level translate to gains in real applications. This section introduces two critical statistical applications and their hardware accelerators used in the case study.

\subsection{Algorithms and Applications}
\label{sec:case_intro}

\textbf{HMM and VICAR}.
The Hidden Markov Model (HMM) is a widely used statistical model~\cite{hmm_tutorial,eddy1996hidden,durbin1998biological,mor2021systematic,
firtina2022aphmm,liu2021variational}. An HMM consists of a sequence of states $q={q_0, q_1, ..., q_{T-1}}$ and a sequence of observations $O={O_0, O_1, ..., O_{T-1}}$ over time length $T$, where $q_t$ and $O_t$ are the hidden state and observation at time $t$, respectively. The transition matrix ($A$) and the emission matrix ($B$) are input probabilities. $A(i,j)$ is the probability from $q_i$ to $q_j$. $B(i,j)$ is the probability of observing $O_j$ given $q_i$. An HMM is denoted as $\lambda=(A, B)$.

The forward algorithm in HMM is used to compute the likelihood $P(O|\lambda)$ of observing $O$ under $\lambda$~\cite{hmm_tutorial}.
Listing~\ref{lst:hmm_lst} shows the algorithm, where $A$ and $B$ probabilities are iteratively multiplied and accumulated. $alpha$ elements are output probabilities that decreases over $t$ as in Figure~\ref{fig:exp_over_time}.

VICAR is a novel phylogenetics tool to analyze evolutionary parameters of species trees using HMM and the forward algorithm~\cite{liu2021variational}. In VICAR, likelihoods are extremely small, as low as $2^{-2,900,000}$ when running on $T=500,000$ site-long HCG sequences. These extremely small numerical values must be preserved because underflow to zero would prevent proper convergence and lead to incorrect results.

\textbf{PBD and LoFreq}. Poisson Binomial Distribution (PBD) is a powerful statistical tool widely used in applications such as finance and bioinformatics~\cite{hong2013computing,schumacher1999binomial,wilm2012lofreq,pbd_useful,tang2023poisson,tejada2011role,rosenman2019some,salas2018techniques}. PBD is a probability distribution describing independent Bernoulli trials. Each trial $i$ has a binary outcome (success or failure) and a prior success probability $p_{i}$. PBD is used to model $N$ trials where $K$ successes are observed.

The key procedure is to compute the probability mass function (PMF) and p-value.
The PMF is defined by $Pr_{n}(X\!\!=\!\!k)$, which is the probability of having exactly $k$ successes in the first $n$ trials. P-values are critical in statistical hypothesis testing and are computed from the PMF. A p-value smaller than a predefined threshold indicates a significant result~\cite{di2020statistical}. The algorithm is shown in Listing~\ref{lst:pbd_lst}, which also iteratively sums the products of probabilities.

\begin{minipage}{0.96\linewidth}
\begin{lstlisting}[language=Python,caption={The forward algorithm.},label={lst:hmm_lst}]
for t from 1 to T: # outer loop
    ot = O[t]
    for q from 0 to H: # inner loop
        path_sum = 0
        for p from 0 to H: # innermost loop
            a_prob = A[p][q]
            term = alpha_prev[p] * a_prob
            path_sum += term

        b_prob = B[q][ot]
        alpha[q] = path_sum * b_prob
    alpha_prev = alpha # copy data to alpha_prev
likelihood = sum(alpha)
return likelihood
\end{lstlisting}
\end{minipage}

\begin{minipage}{0.96\linewidth}
\begin{lstlisting}[language=Python,caption={Computing PMF and p-value.},label={lst:pbd_lst}]
for n from 1 to N: # outer loop
    pn = success_probs[n]
    for k from 1 to K: # inner loop
        pr[k] = pr_prev[k] * (1-pn) + pr_prev[k-1] * pn
    pr[0] = pr_prev[0] * (1-pn)
    if n > K:
        pvalue = pvalue_prev + pr_prev[K-1] * pn
    pr_prev = pr # copy data to pr_prev
    pvalue_prev = pvalue
return pvalue
\end{lstlisting}
\end{minipage}

LoFreq is a critical genomics tool that identifies genome variants using PBD to analyze genome alignment data (\textit{columns})~\cite{wilm2012lofreq}. Each column is modeled with PBD, from which a p-value is computed. Each column has its own $N$, $K$, and $success\_probs$. LoFreq determines that genome variants exist in a column if its p-value is $<2^{-200}$. We analyzed the p-values of $222,131$ columns in SARS-CoV-2 data from~\cite{lofreq_trets}. Their p-values span a wide range. Among all, $16,205$ columns are reported to have genome variants (critical). $40\%$ and $5\%$ of these critical columns have a p-value $<2^{-1,074}$ and $2^{-10,000}$, respectively. The smallest observed p-value is $2^{-434,916}$.

The forward algorithm and PBDs are commonly implemented using logarithms~\cite{hmm_numerical,lofreq_pbd_impl,stan_pbd_impl,liu2021variational,durbin1998biological}. Otherwise, these critical likelihoods and p-values will underflow, leading to catastrophic results. Listing~\ref{lst:hmm_lst_log} shows the forward algorithm in log-space, where $LSE$ implements Equation~\eqref{eq:lse}. A non-accumulative add becomes a binary LSE as in Equation~\eqref{eq:lse_two}. $ln\_A$ and $ln\_B$ are pre-computed logarithms of $A$ and $B$.

\subsection{Accelerator Design and Implementation}
\label{sec:case_acc}
The computation in both applications has common characteristics and is amenable to hardware acceleration. For LoFreq, we use the \textit{column unit}, an FPGA accelerator implemented in~\cite{lofreq_trets}. It is able to deliver up to $51.7\times$ end-to-end speedup, which is the fastest accelerator for LoFreq as far as we know. For HMM and VICAR, we implemented an FPGA accelerator, referred to as \textit{forward algorithm unit}. Both accelerators implement computations in log-space. Both are highly parallelized. Compared to the C implementation on the CPU, the forward algorithm unit achieves $66\times$ and $115\times$ speedup when $H=64$ and $128$, respectively. Meanwhile, these FPGA accelerators produce bit-equivalent results to the original CPU software.

\begin{minipage}{0.96\linewidth}
\begin{lstlisting}[language=Python,caption={The forward algorithm using logarithms.},label={lst:hmm_lst_log}]
for t from 1 to T: # outer loop
    ot = O[t]
    for q from 0 to H: # inner loop
        terms = []
        for p from 0 to H: # innermost loop
            a_prob = ln_A[p][q]
            term = alpha_prev[p] + a_prob
            terms[p] = term
        path_sum = LSE(terms)
        b_prob = ln_B[q][ot]
        alpha[q] = path_sum + b_prob
    alpha_prev = alpha # copy data to alpha_prev
log_likelihood = LSE(alpha)
return log_likelihood
\end{lstlisting}
\end{minipage}

The first characteristic is a nested loop structure where the outer loop is sequential but the inner loop is parallel. As shown in Section~\ref{sec:case_intro}, the outer loop is sequential due to data dependency: $alpha$ and $pr$ are iteratively computed from $alpha\_prev$ and $pr\_prev$, respectively. In contrast, the inner loop iterations are independent from each other and thus can be parallelized. One notable challenge in the forward algorithm is the accumulation in the innermost loop (line 8 in Listing~\ref{lst:hmm_lst}), which prevents straightforward parallelization.

\begin{figure}[t]
    \centering
    \subfigure[Logarithm-based PE.]
    {
        \centering
        \includegraphics[width=1.0\columnwidth]{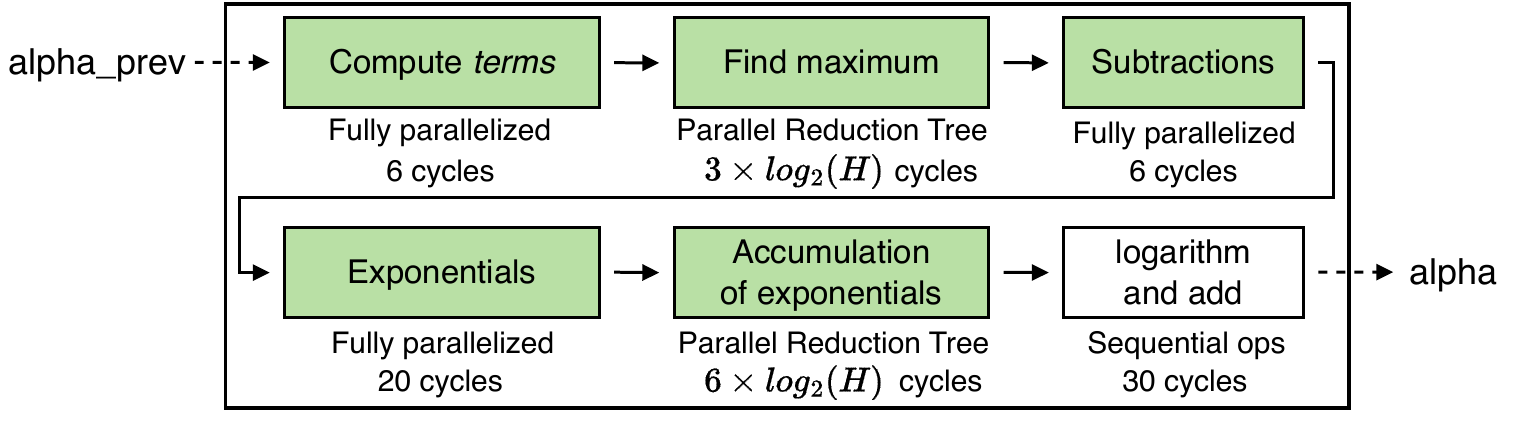}
        \label{fig:hmm_pe_log}
    }
    \subfigure[Posit-based PE.]
    {
        \centering
        \includegraphics[width=1.0\columnwidth]{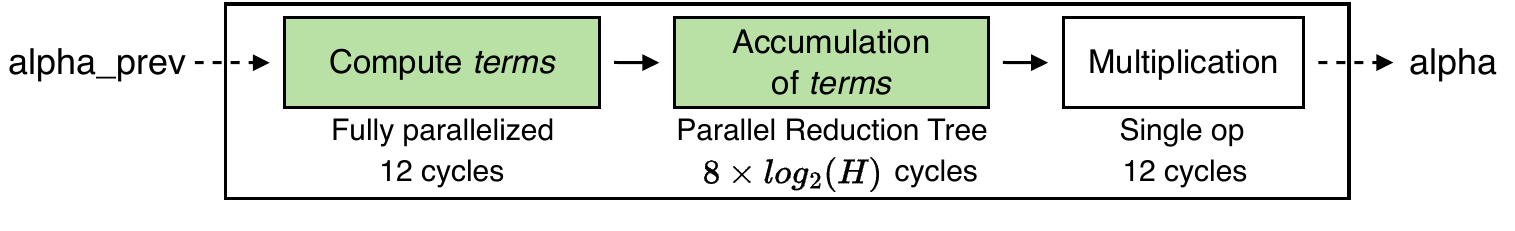}
        \label{fig:hmm_pe_posit}
    }
    \caption{Processing Element (PE) in forward algorithm units.}
    \label{fig:hmm_pe_design}
\end{figure}

Both accelerators implement the computation of an inner loop iteration in Processing Elements (PE). PEs in both accelerators are fully pipelined, initiating a new \textit{inner loop} iteration every clock cycle. Note that the PE in the forward algorithm unit fully parallelizes the \textit{innermost loop}, as shown in Figure~\ref{fig:hmm_pe_log}, in order to be fully pipelined. Thus, these PEs are hardwired for fixed $H$ parameters.

The second characteristic is that the sequence input is long: $O$ is of length $T$ and $success\_probs$ is of length $N$. Both $T$ and $N$ are large, making the input hard to fit in on-chip SRAMs. In every outer loop iteration, a new sequence element is accessed for one use. Thus, both accelerators store the input in DRAM and implement a prefetcher for DRAM accesses.

\begin{figure}[t]
    \centering
    \includegraphics[width=1.0\columnwidth]{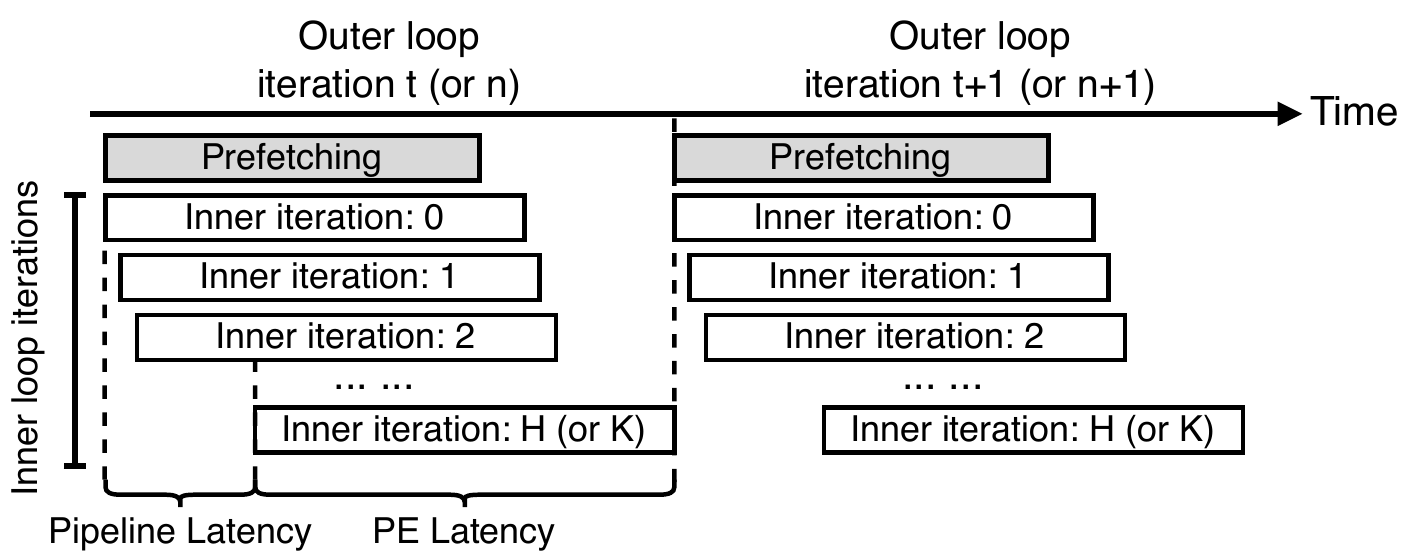}
    \caption{Execution timeline of the accelerators.}
    \label{fig:acc_exec}
    \end{figure}

In both accelerators, the fully pipelined PEs and the prefetcher run in parallel. Figure~\ref{fig:acc_exec} visualizes their execution. The total number of cycles of both is computed by $\textit{outer loop bound} \times (\textit{pipeline latency} + \textit{PE latency})$. The $\textit{outer loop bound}$ is $T$ and $N$ for VICAR and LoFreq, respectively.  The pipeline latency is the number of cycles in initiating new \textit{inner loop} iterations ($H$ for VICAR and $K$ for LoFreq).

\subsection{Posit-based Implementation}
\label{sec:case_posits}
Using posit in the accelerator reduces the critical path latency. This is clearly visualized in Figure~\ref{fig:hmm_pe_log}. In a forward algorithm unit, a log-based PE has to implement an $H$-nary LSE unit which contains $H$ exponential units, $H$ adders, $H/2$ comparators, and one logarithm unit. The PE function stages are shown in Figure~\ref{fig:hmm_pe_log}, and the latency is $62 + 9 \times log_2(H)$ cycles. In contrast, when using posits, the PE's logic is largely simplified. Its latency becomes $24 + 8 \times log_2(H)$ cycles, with a reduction of $38 + log_2(H)$ cycles, as shown in Figure~\ref{fig:hmm_pe_posit}. A log-based PE implements an adder and a binary LSE unit. Its latency is 73 cycles: 64 cycles for an LSE,  6 cycles for an add, and 3 cycles for conditional logic. In a posit-based PE, the latency is reduced to $30$ cycles.

Besides, the posit-based accelerators consume less than half of the resources used by their logarithm-based counterparts. Such a resource-saving enables building more hardware units on the FPGA, leading to even greater speedup. Moreover, using posit shifts the performance bottleneck from the PEs to the prefetcher when $H$ (or $K$) is small, unlocking opportunities for further speedup by reducing DRAM access latency.

\section{Evaluation}
\label{sec:eval}

\subsection{Evaluation Methodology}
\label{sec:setup}
\textbf{System Setup}. The evaluation focuses on hardware accelerators implemented on an Xilinx Alveo U250 card. All designs were developed in HLS C++. Arithmetic operations in the logarithm-based accelerators are implemented using the standard Xilinx library (LogiCORE IP v7.1~\cite{logicoreip}). Arithmetic operations in the posit-based accelerators are implemented using the MArTo HLS C++ library~\cite{uguen2019evaluating}. Xilinx Vitis 2020.02 was used for synthesis, placement, and routing.

\textbf{Metrics and Baselines}. The posit-based and logarithm-based accelerators are compared in \textit{performance}, \textit{resource cost}, and \textit{numerical accuracy}. The logarithm-based accelerators are the baselines. \textbf{Recall from section~\ref{sec:case_acc} that these baselines are highly optimized FPGA accelerators.} Each column unit has $8$ PEs, and each forward algorithm unit has one PE that is fully pipelined and completely parallelizes the innermost loop. The wall clock execution time of accelerators is measured for performance comparison. To minimize the impact of implementation differences between MArTo and the Xilinx IP, all accelerators are implemented to operate at 300 MHz for evaluation.  Meanwhile, the maximum achievable clock frequency is also shown in Table~\ref{tab:fau_total_resource} and~\ref{tab:total_resource} for reference.

The numerical accuracy of final application-level results is evaluated. Baseline correct results are computed using 256-bit MPFR~\cite{mpfr_paper}. Relative errors between the accelerator results and the correct results are used to measure accuracy.

\textbf{Datasets}. To evaluate the forward algorithm units, both Human-Chimp-Gorilla (HCG) data~\cite{hcg_data,liu2021variational} and synthetic HMM data are used. For the HCG data, input probabilities $A$ and $B$ are generated by VICAR. For the synthetic HMM data, $A$ and $B$ are synthesized from the Dirichlet distribution, and $O$ is universally sampled. Eight typical SARS-CoV-2 datasets (average of $N$ is $309,189$) from~\cite{lofreq_trets} are used to evaluate the column units. In these eight datasets, there are in total $222,131$ columns, among which $16,205$ are critical (p-value $< 2 ^ {-200}$). In these datasets, $N$ and $K$ are diversely distributed, unlike $T$ and $H$ in VICAR. These real datasets are used for both performance and numerical accuracy evaluation.

\begin{figure}[t]
  \subfigure[Wall clock execution time ($T=500,000$).]
  {
      \centering
      \includegraphics[width=1.0\columnwidth]{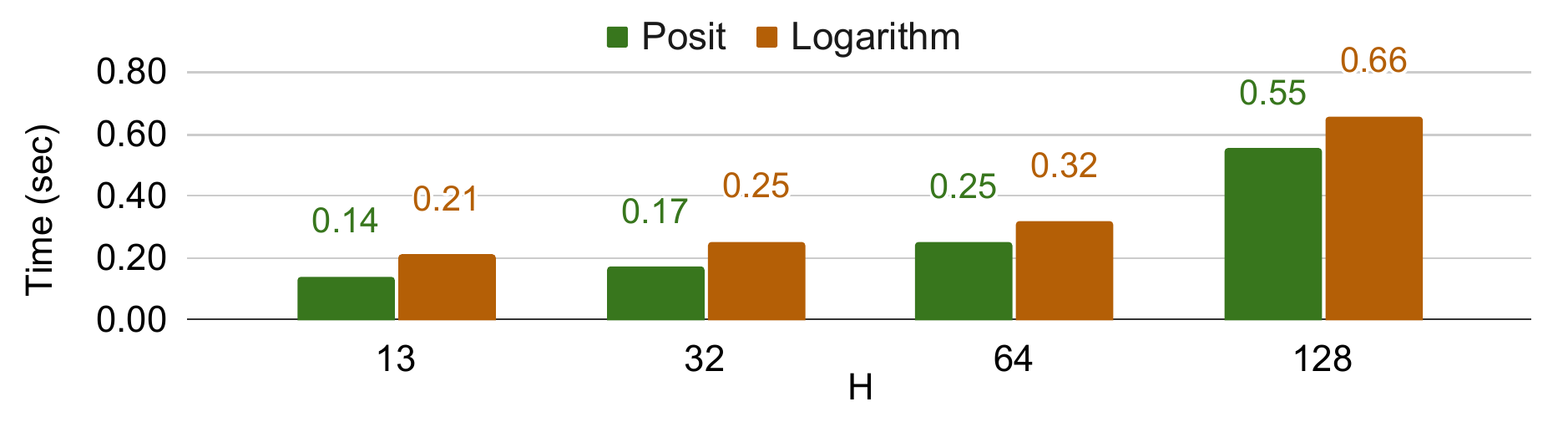}
      \label{fig:fau_time}
  }
  \subfigure[Relative improvement.]
  {
      \centering
      \includegraphics[width=1.0\columnwidth]{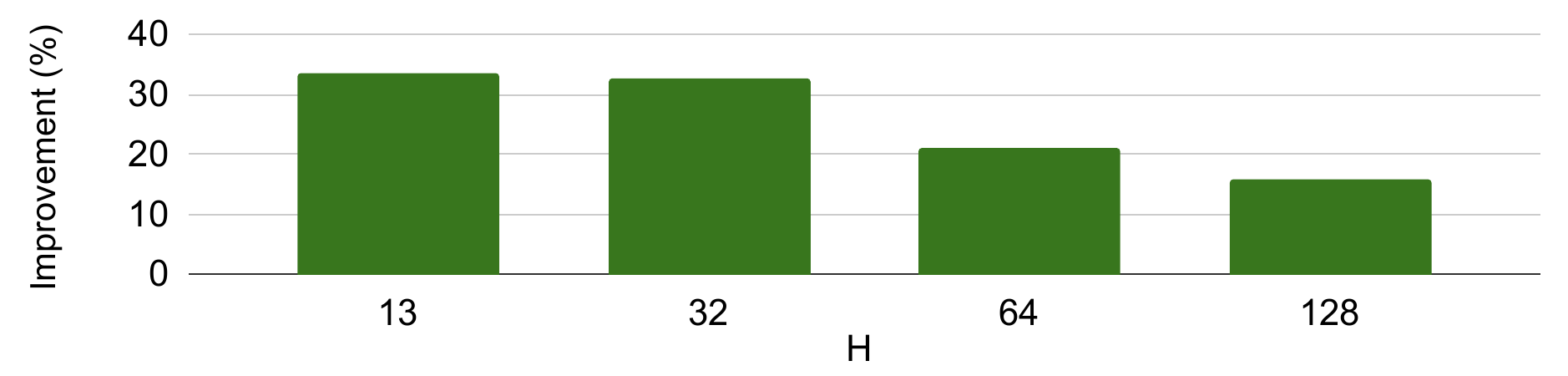}
      \label{fig:fau_ratio}
  }
  \caption{Performance of forward algorithm units.}
  \label{fig:hmm_performance}
\end{figure}

  \begin{figure}[t]
    \subfigure[Wall clock execution time.]
    {
        \centering
        \includegraphics[width=1.0\columnwidth]{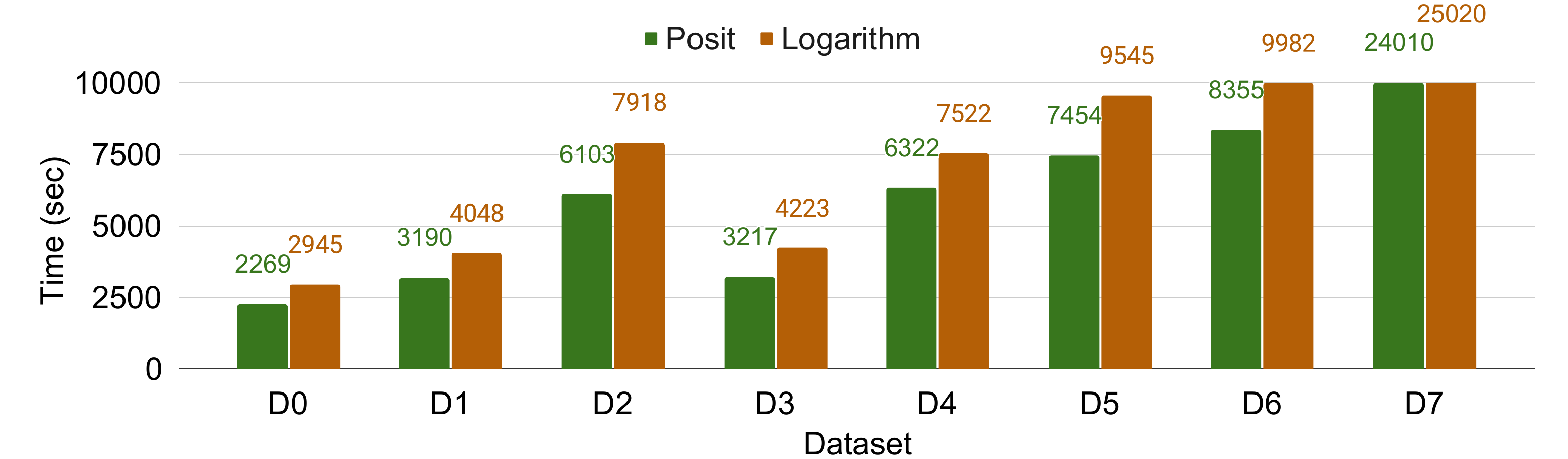}
        \label{fig:cu_time}
    }
    \subfigure[Relative improvement.]
    {
        \centering
        \includegraphics[width=1.0\columnwidth]{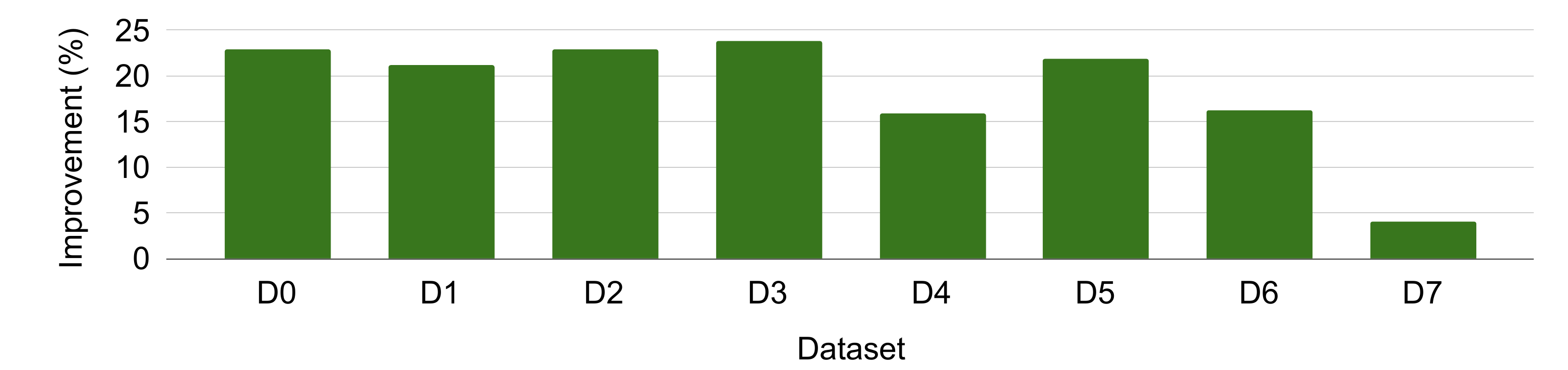}
        \label{fig:cu_ratio}
    }
    \caption{Performance of column units.}
    \label{fig:lofreq_performance}
\end{figure}

\begin{center}
\begin{table}[t]
  \centering
  \caption{Resource Use of Forward Algorithm Units.}
  \label{tab:fau_total_resource}
  \resizebox{1.0\columnwidth}{!}
  {
  \begin{tabular}{r | r | r r r r r | r}
    \toprule
      & H & CLB & LUT & Register & DSP & SRAM & \vtop{\hbox{\strut Max Clock }\hbox{\strut Frequency (MHz)}}\\ [0.5ex]
    \midrule
    Logarithm  & 13 & 14,308 & 68,966 & 61,720 & 275 & 43 & 345 \\
    posit(64,18) & 13 & 6,272 & 26,093 & 32,271 & 143 & 43 & 330 \\
    \midrule
    Reduction &  & 56.16\% & 62.16\% & 47.71\% & 48.00\% & 0 & 4.35\% \\
    \midrule
    \midrule
    Logarithm  & 32 & 27,264 & 145,300 & 119,435 & 560 & 98 & 345 \\
    posit(64,18) & 32 & 12,090 & 55,910 & 67,906 & 314 & 102 & 330 \\
    \midrule
    Reduction &   & 55.66\% & 61.52\% & 43.14\% & 43.93\% & -4.08\% & 4.35\% \\
    \midrule
    \midrule
    Logarithm  & 64 & 47,058 & 273,525 & 216,083 & 1,021 & 250 & 332 \\
    posit(64,18) & 64 & 23,187 & 103,948 & 125,875 & 602 & 258 & 330 \\
    \midrule
    Reduction &  & 50.73\% & 62.00\% & 41.75\% & 41.04\% & -3.20\% & 0.61\% \\
    \midrule
    \midrule
    Logarithm  & 128 & 50,690 & 308,719 & 258,834 & 1,040 & 1,406 & 308 \\
    posit(64,18) & 128 & 23,775 & 123,011 & 157,696 & 602 & 1,410 & 300 \\
    \midrule
    Reduction &  & 53.10\% & 60.15\% & 39.07\% & 42.12\% & -0.28\% & 2.67\% \\
    \bottomrule
  \end{tabular}
  }
\end{table}
\end{center}

\begin{center}
  \begin{table}[t]
    \centering
    \caption{Resource Use of Column Units.}
    \label{tab:total_resource}
    \resizebox{1.0\columnwidth}{!}
    {
    \begin{tabular}{r | r | r r r r r | r }
      \toprule
        & \# of PEs & CLB & LUT & Register  & DSP & SRAM & \vtop{\hbox{\strut Max Clock }\hbox{\strut Frequency (MHz)}}\\ [0.5ex]
      \midrule
      Logarithm  & 8 & 15,476 & 75,894 & 76,300 & 386 & 236 & 341 \\
      posit(64,12) & 8 & 8,619 & 27,270 & 37,963 & 153 & 258 & 330 \\
      \midrule
      Reduction & - & 44.31\% & 64.07\% & 50.25\% & 60.36\% & -9.32\% & 3.22\% \\
      \bottomrule
    \end{tabular}
    }
  \end{table}
  \end{center}

\begin{figure}[t]
  \centering  \includegraphics[width=1.0\columnwidth]{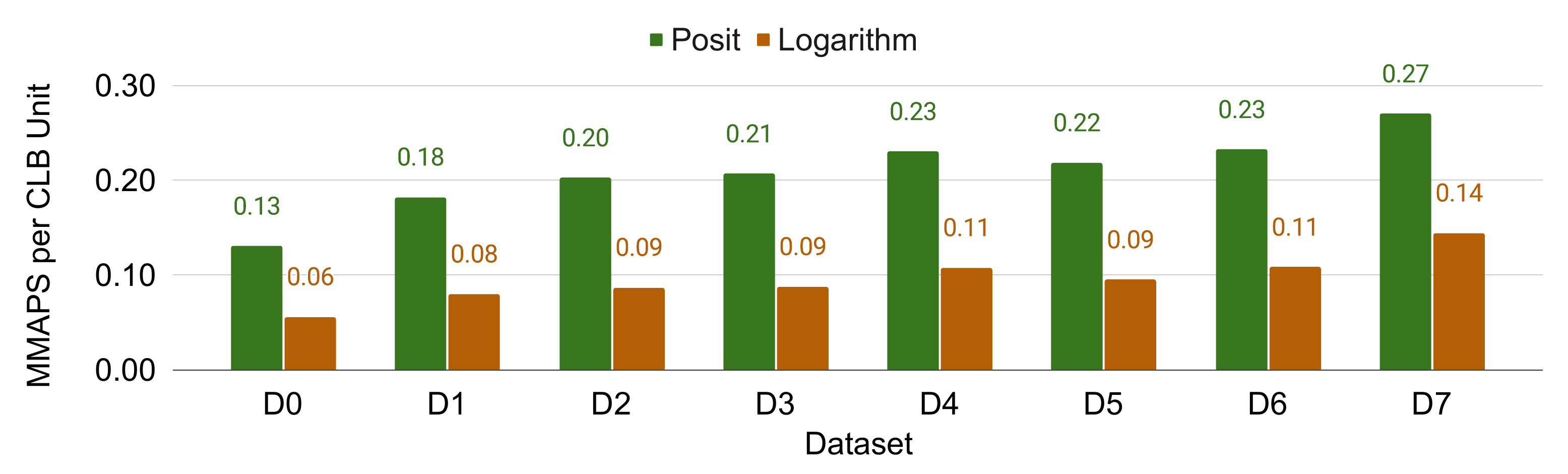}
  \caption{Performance per Resource Unit.}
  \label{fig:perf_per_resource}
\end{figure}

\begin{figure*}[t]
  \centering
  \includegraphics[width=2.0\columnwidth]{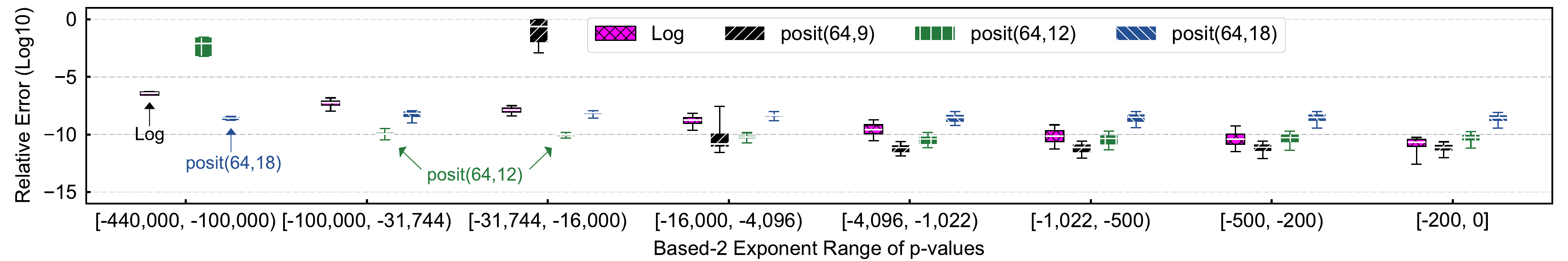}
  \caption{Accuracy of final p-values of different magnitudes.}
  \label{fig:pvalue_boxplot}
\end{figure*}

\begin{figure}[t]
\centering
\subfigure[$T=100,000$.]
{
    \centering
    \includegraphics[width=0.45\columnwidth]{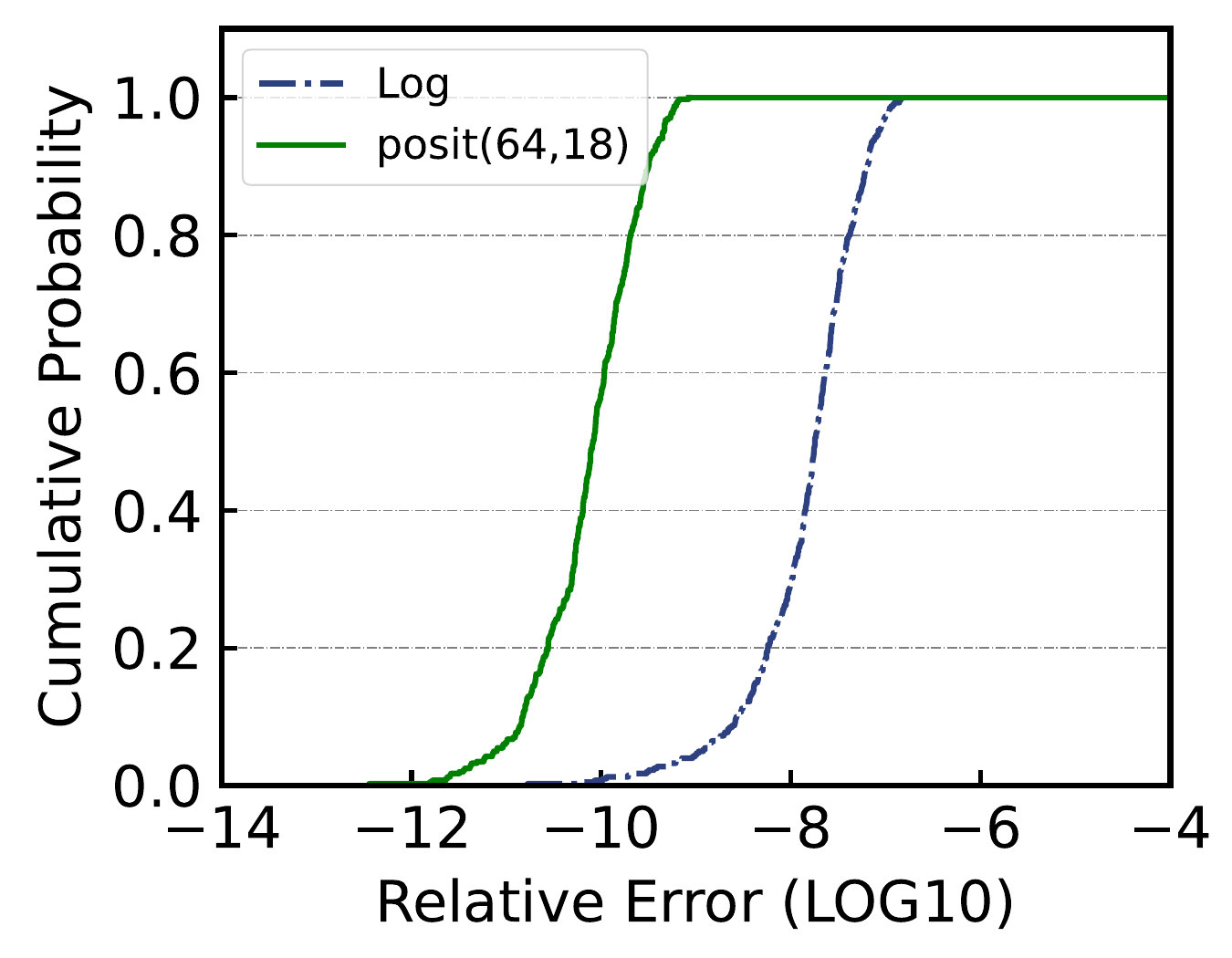}
    \label{fig:likelihood_100000}
}
\subfigure[$T=500,000$.]
{
    \centering
    \includegraphics[width=0.45\columnwidth]{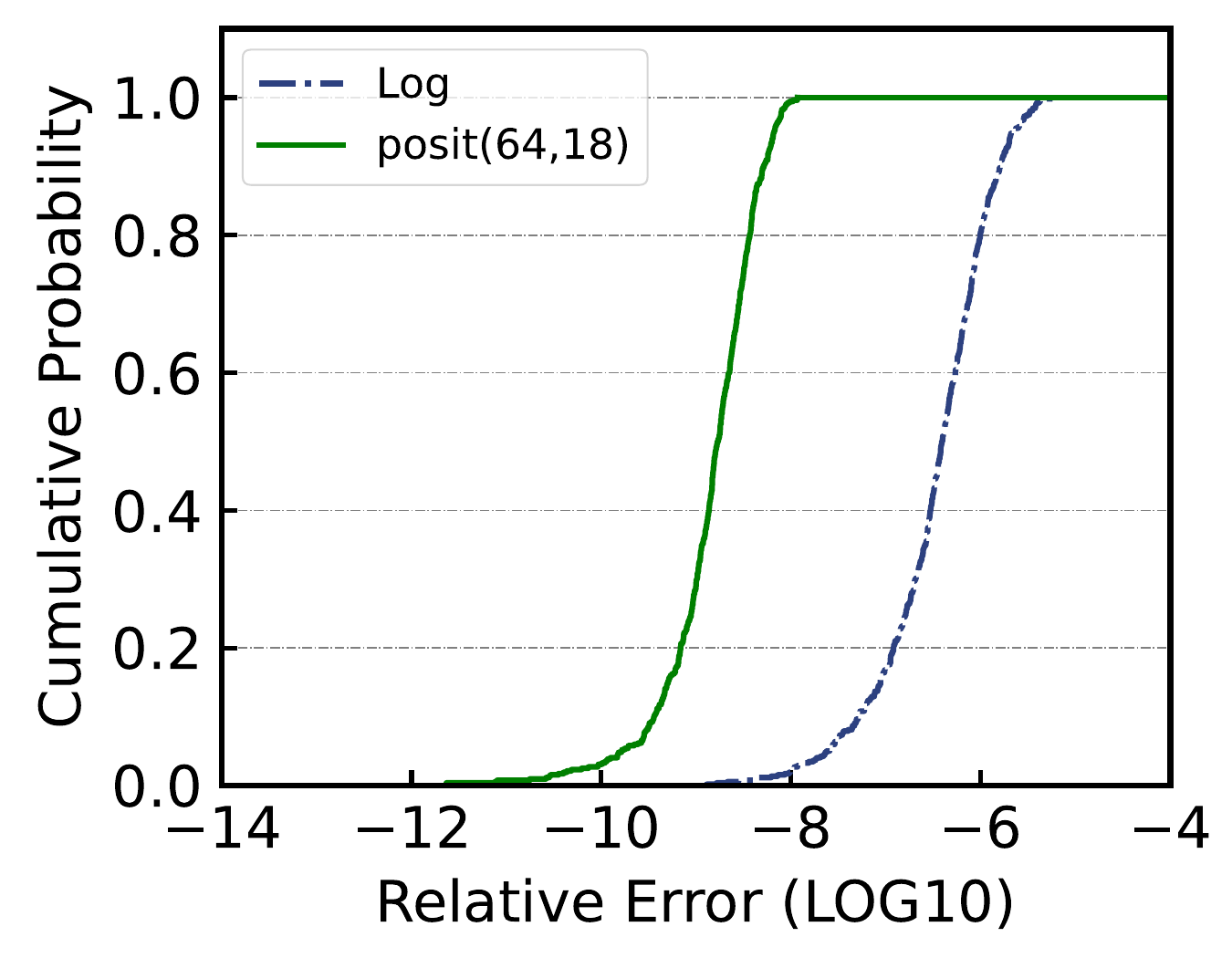}
    \label{fig:likelihood_500000}
}
\caption{Overall accuracy of final likelihoods in VICAR.}
\label{fig:likelihood_cdf}
\end{figure}

\begin{figure}[t]
  \centering
  \subfigure[p-values $< 2 ^ {-200}$ (critical).]
  {
      \centering
      \includegraphics[width=0.45\columnwidth]{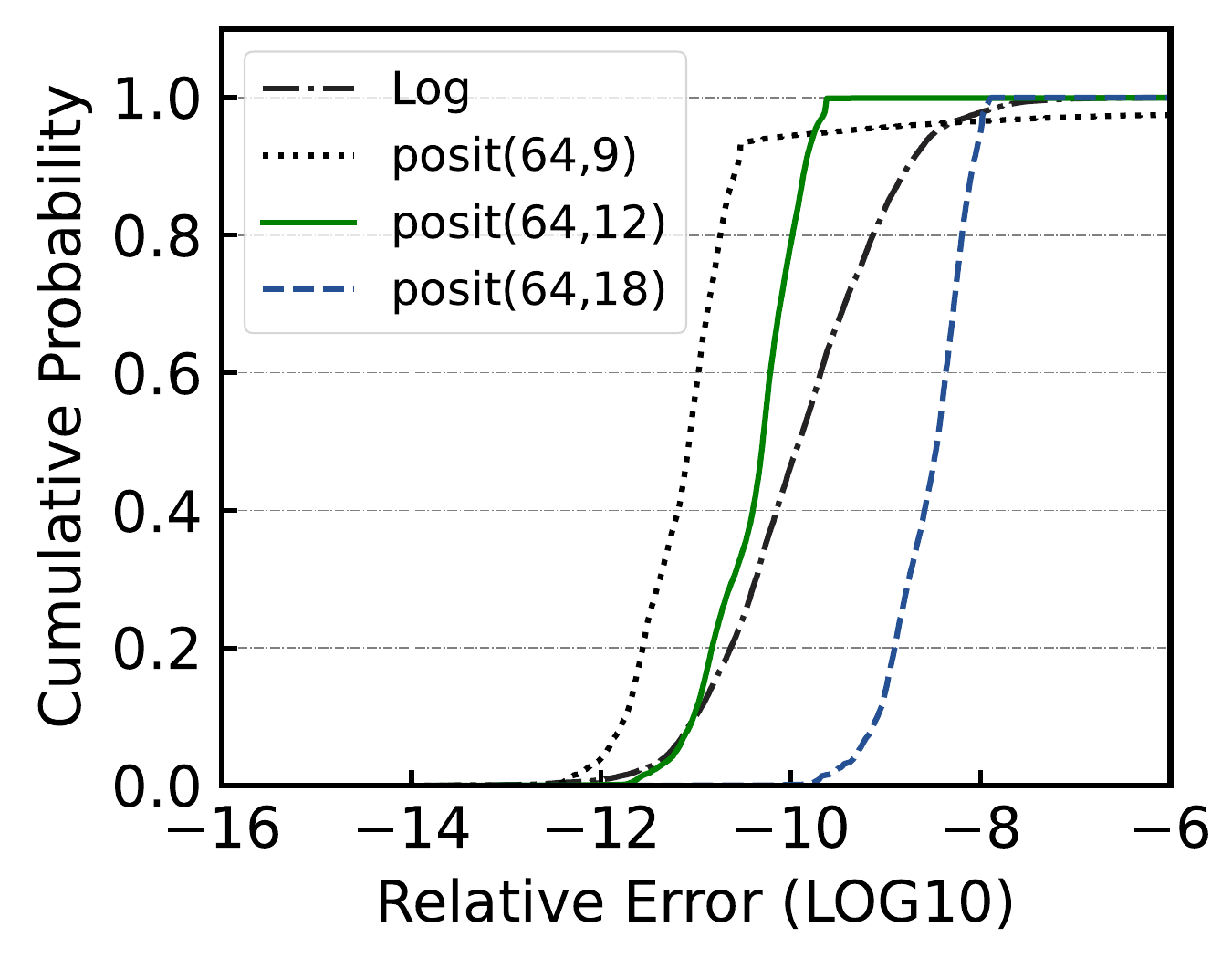}
      \label{fig:pvalue_snvs}
  }
  \subfigure[p-values $>= 2 ^ {-200}$ (less important).]
  {
      \centering
      \includegraphics[width=0.45\columnwidth]{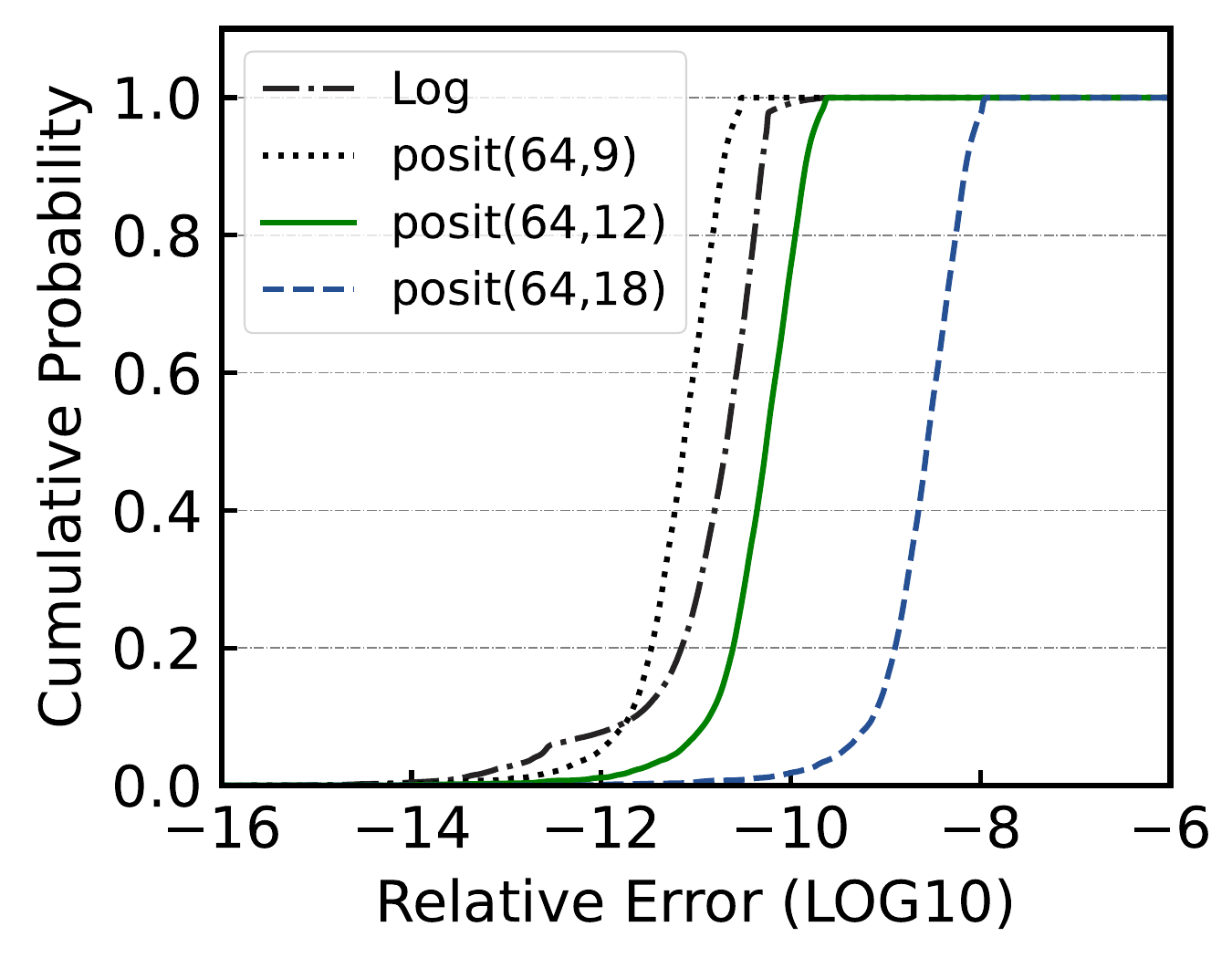}
      \label{fig:pvalue_large}
  }
  \caption{Overall accuracy of final p-values in LoFreq.}
  \label{fig:pvalue_accuracy_cdf}
\end{figure}

\subsection{Single Unit Performance}
\textbf{A single posit-based unit is consistently $15\%$ to $33\%$ faster than its logarithm-based counterpart.} Figures~\ref{fig:hmm_performance} and~\ref{fig:lofreq_performance} highlight the performance of posit-based and log-based units. The relative improvement is calculated by the execution time reduction divided by the logarithm's execution time. The relative speedup tends to be small when $H$ or $K$ is large. This is because the improvement from posits is small relative to the large pipeline latency, as described in Section~\ref{sec:case_acc}.

\subsection{Performance Per Resource Unit}
\textbf{Using posits leads to $60\%$ lower resource use.} Tables~\ref{tab:fau_total_resource} and~\ref{tab:total_resource} show the resource cost of forward algorithm units and column units, respectively. CLB (Configurable Logic Block) is the resource building block on Xilinx UltraScale+ FPGAs~\cite{clb_doc}. Each CLB slice contains both LUTs and registers. Posit-based accelerators consume only about $40\%$ of LUTs and $60\%$ of registers and DSPs as used by their log-based counterparts.

\textbf{Such resource reduction not only yields an area advantage, but also translates to a larger speedup}, when more units are implemented, enabling more parallelization. For example, an FPGA die slice (SLR) on a U250 can implement at most 4 log-based column units. In contrast, it can easily fit 10 posit-based column units.

\textbf{Using posits leads to $2\times$ higher performance per resource unit compared to using logarithms.} The column units are used for analysis. The performance metric is the throughput of the multiply-and-add operations shown in line 4 of Listing~\ref{lst:pbd_lst}, which is dubbed \textit{MMAPS}, short for \textit{Million Multiplies and Add Per Second}. Each dataset used has about $10^{13}$ multiply-and-add operations. The resource metric is the CLB usage as shown in Table~\ref{tab:total_resource}. Overall, \textit{MMAPS per CLB Unit} is the metric to measure performance per resource unit. As shown in Figure~\ref{fig:perf_per_resource}, posit-based column units perform \textbf{twice} as many MMAPS per CLB unit on all datasets.

\subsection{Application Numerical Accuracy}
\textbf{Using posits results in higher accuracy than using logarithms.} The higher accuracy at the arithmetic level, as described in Section~\ref{sec:tradeoff_acc}, translates to more accurate application-level final results. However, posits still suffer from poor accuracy when numbers are close to or outside the range, since their ranges are still limited.

\textbf{VICAR.} The CDFs in Figure~\ref{fig:likelihood_cdf} present the overall accuracy of the final likelihoods computed by the accelerators. When $T=100,000$ and $T=500,000$, the likelihoods are approximately $2^{-590,000}$ and $2^{-2,900,000}$, respectively. Each CDF shows the accuracy distribution of the final likelihoods using $512$ different $A$ and $B$ input matrices ($128$ for each $H$). The overall accuracy is higher when the curve is more skewed towards the left. For example, Figure~\ref{fig:likelihood_500000} shows that $100\%$ posit$(64,18)$ results have a relative error $<10^{-8}$, compared to only $2.4\%$ logarithm results achieving the same. On such extremely small numbers, using posit$(64,18)$ leads to two orders of magnitude higher accuracy compared to logarithms.

\textbf{LoFreq.} Unlike the VICAR likelihoods, LoFreq p-values span an extremely wide range from $2^{-434,916}$ to $1.0$. Similar to Figure~\ref{fig:likelihood_cdf}, CDFs in Figure~\ref{fig:pvalue_accuracy_cdf} show the overall accuracy across all p-values. Figure~\ref{fig:pvalue_boxplot} shows the accuracy of application-level final results (p-values) in different magnitudes. Note that extreme cases with relative error $>=1$ are \textbf{not} included in Figure~\ref{fig:pvalue_boxplot}. This is why posit$(64,9)$ is absent in the two leftmost ranges of the figure, and binary64 is not shown at all.

For most p-values, posit$(64,12)$ and posit$(64,9)$ have higher accuracy than using logarithms. In Figure~\ref{fig:pvalue_snvs}, curves of posit$(64,9)$ and posit$(64,12)$ lie to the left of the logarithm's curve, indicating higher overall accuracy. In particular, $99\%$ posit$(64,12)$ results have a relative error $<10^{-10}$, while only $60\%$ logarithm results achieve such accuracy. In Figure~\ref{fig:pvalue_large}, posit$(64,9)$ achieves the highest accuracy on the non-critical p-values. \textbf{This shows how the arithmetic accuracy gains shown in Section~\ref{sec:tradeoff_acc} translate directly to more accurate application-level final results.}

However, posits do not always achieve higher accuracy due to limited ranges. Underflow still occurs: among all $222,131$ p-values, posit$(64,9)$ and posit$(64,12)$ underflow on $132$ and $2$ of them, respectively. Posit$(64,18)$ does not underflow. Second, relative errors can exceed $1$ when values approach the format's representable limit. This is because most bits are used as regime, not fraction. posit$(64,9)$ and posit$(64,12)$ have $30$ and $2$ such high-error cases, respectively. The largest observed relative error of posit$(64,9)$ and posit$(64,12)$ is about $10^{295}$ and $10^{2,129}$, respectively. In contrast, posit$(64,18)$ has zero such cases thanks to its large \textit{useed} and wide range.


The trade-off among the three posits is shown in Figure~\ref{fig:pvalue_boxplot}. The accuracy of posit$(64,9)$ is the highest on p-values $> 2^{-16,000}$ but drastically drops and eventually underflows on smaller values. Posit$(64,12)$ achieves high accuracy on all p-values $> 2^{-100,000}$, much wider than the high accuracy range of posit$(64,9)$. However, posit$(64,12)$ still underflows on extreme cases like $2^{-434,916}$. In contrast, while posit$(64,18)$ has the worst accuracy on p-values $> 2^{-16,000}$, it has the highest accuracy on the extremely small p-values, such as $2^{-434,916}$. On numbers of that magnitude, posit$(64,18)$'s accuracy is notably higher than using logarithms.

\section{Related Works}
\label{sec:related}

\textbf{Posits}.
The performance and resource cost of  posit arithmetic on FPGAs have been studied in~\cite{jaiswal2019pacogen,podobas2018hardware,
uguen2019evaluating,forget2021comparing,chaurasiya2018parameterized,8342187,
murillo2022comparing,murillo2020customized,electronics9101622,de2019posits,
shekhawat2023phac,lehoczky2018high,johnson2018rethinking,murillo2023plaus,
tiwari2021peri,mallasen2022percival,mallasen2023big,sharma2023clarinet,mallasen2022customizing}. Some prior works have compared the accuracy of posits with IEEE floating-point numbers~\cite{debug_posit_pldi,omtzigt2022universal, tiwari2021peri,
chien2019posit,mallasen2022percival,mallasen2023big,murillo2023generating,
van2019accelerator,klowertowards,murillo2023plaus,
lim2021high,lim2021approach}. A few prior works have built accelerators using posits~\cite{sommer2020comparison,van2018enabling,murillo2023generating}.

Our paper stands out as the first to comprehensively study numerical accuracy comparing posits, binary64, and, particularly, using logarithms, and the first to propose using posits in statistical computations operating on small numbers.

\textbf{Statistical Accelerators}.
Prior works have shown promising results in accelerating statistics~\cite{
chai2022coopmc,ko20203mm,ko2020scalable,ko2019flexgibbs,
banerjee2019acmc,zhang2021statistical,zhang2020accelerating,fan2021high} and bioinformatics~\cite{
huang2017hardware,Salamat2020FPGAAO,subramaniyan2021accelerated,
alser2020sneakysnake,guo2019hardware,lo2020algorithm,
gu2023gendp,jiang2021exma,turakhia2018darwin}. In contrast, our paper
targets a class of statistical computations that have not been studied, and is the first to propose using posits in such computations.

\textbf{Other Formats}. Logarithmic Number System (LNS) encodes log values in fixed-point rather than floating-point~\cite{alam2021low,parhami2020computing}. It was designed for low-precision, narrow-range applications~\cite{parhami2020computing} and low-precision training~\cite{zhao2022lns,haghi2024bridging}, which typically use 16 or fewer bits. For these applications, LSE in LNS can be efficiently implemented using lookup tables with pre-computed $log(1+exp(x))$, eliminating expensive logarithms and exponentials.

In contrast, LNS is not suitable for wide-range, high-accuracy statistical computations. This is because the lookup table optimizations are impractical for 64-bit numbers.

Rescaling is another approach that prevents underflow by multiplying small numbers with a constant scaling factor~\cite{li2015robust,murphy2002hidden}. However, it is impractical in our target applications, where numbers span an extremely wide range.


\section{Conclusion}
\label{sec:conc}

Due to the range limitations of IEEE floating-point numbers, statistical computations are often done in log-space. This paper quantitatively reveals that using logarithms not only leads to worse performance and resource cost, but also harms numerical accuracy. Furthermore, this paper has built FPGA accelerators to show that using posits achieves better accuracy, lower resource cost, and higher performance. The key insight is that posits succeed by allowing bits to be dynamically allocated between exponent and fraction bits, instead of using the fixed allocation found in IEEE floating-point. In particular, using posits makes the final application-level results two orders of magnitude more accurate. Compared to log-based accelerators, posit-based accelerators achieve up to 60\% lower resource use and up to 33\% higher performance. This results in gains of up to a factor of 2 in performance per resource unit on the FPGA.

\bibliographystyle{IEEEtranS}
\bibliography{reference}

\end{document}